\documentclass[12pt]{amsart}
\usepackage{amsmath, amssymb, latexsym, amsthm, bm}
\usepackage{mathrsfs}
\usepackage[all]{xy}

\newcommand{\Pa}[9]{\bibitem{#1} {#2}, \emph{#3}, {#4} \textbf{#5} ({#6}), {#7}--{#8}.}
\newcommand{\ed}{

\end{document}
}

\newcommand{\alephes}{{\aleph_0}}

\newcommand{\arx}[1]{\texttt{arxiv.org/math/#1}}
\newcommand{\bq}{\begin{quote}}
\newcommand{\eq}{\end{quote}}

\newcommand{\inv}{^{-1}}

\newcommand{\N}{\mathbb{N}}

\newcommand{\seq}[1]{\{#1\}_{n\in\N}}
\newcommand{\sseq}[1]{\setseq{#1}}
\newcommand{\setseq}[1]{\{#1 : n\in\N\}}

\newcommand{\scrA}{\mathscr{A}}
\newcommand{\scrB}{\mathscr{B}}

\newcommand{\cB}{\mathcal{B}}

\newcommand{\cF}{\mathcal{F}}

\newcommand{\cO}{\mathcal{O}}

\newcommand{\cU}{\mathcal{U}}
\newcommand{\Union}{\bigcup}
\newcommand{\cV}{\mathcal{V}}
\newcommand{\cW}{\mathcal{W}}

\newcommand{\Impl}{\Rightarrow}
\long\def\forget#1\forgotten{}

\newcommand{\oo}{\infty}

\newcommand{\w}{\omega}

\newcommand{\nin}{\not\in}

\newcommand{\sbst}{\subseteq}

\newcommand{\sm}{\setminus}

\newtheorem{thm}{Theorem}

\newtheorem{lem}[thm]{Lemma}
\newtheorem{cor}[thm]{Corollary}

\theoremstyle{definition}

\theoremstyle{remark}

\newcommand{\be}{\begin{enumerate}}
\newcommand{\ee}{\end{enumerate}}
\newcommand{\bi}{\begin{itemize}}
\newcommand{\itm}{\item}
\newcommand{\ei}{\end{itemize}}
\newcommand{\Arh}{Arhangel'ski\u{\i}}

\forget
\forgotten

\newcommand{\sone}{\mathsf{S}_1}
\newcommand{\sfin}{\mathsf{S}_\mathrm{fin}}

\newcommand{\scof}{\mathsf{S}_\mathrm{cof}}
\newcommand{\sinf}{\mathsf{S}_{\oo}}

\newcommand{\capinf}{\mbox{$\bigcap\nolimits_{\oo}$}}

\author{Boaz Tsaban}
\thanks{Supported by the Koshland Center for Basic Research.}
\address{Department of Mathematics, Bar-Ilan University,
Ramat-Gan 52900, Israel;
and
Department of Mathematics,
Weizmann Institute of Science, Rehovot 76100, Israel.
}
\email{tsaban@math.biu.ac.il}
\urladdr{http://www.cs.biu.ac.il/\~{}tsaban}

\subjclass[2000]{54D20, 
54H11. 
}

\keywords{\Arh{} $\alpha_i$ spaces, Ko\v{c}inac $\alpha_i$ selection principles,
totally bounded groups.}

\title{On the Ko\v{c}inac $\alpha_i$ properties}

\begin{document}

\begin{abstract}
The Ko\v{c}inac $\alpha_i$ properties, $i=1,2,3,4$, are
generalizations of \Arh's $\alpha_i$ local properties.
We give a complete classification of these properties when
applied to the standard families of open covers of topological
spaces or to the standard families of open covers of topological
groups. One of the latter properties characterizes
totally bounded groups.
We also answer a question of Ko\v{c}inac.
\end{abstract}

\maketitle

\section{Introduction}

We say that $\cU$ is a \emph{cover} of a set $X$ if $X\nin\cU$
and $X=\Union\cU$.
For topological spaces, various special families of covers have been
extensively studied in the literature, in a framework called \emph{selection
principles}, see the surveys \cite{LecceSurvey, KocSurv, ict}.

The main types of covers are defined as follows.
Let $\cU$ be a cover of $X$. $\cU$ is an \emph{$\omega$-cover} of $X$ if
each finite $F\sbst X$ is contained in some $U\in\cU$.
$\cU$ is a \emph{$\gamma$-cover} of $X$ if
$\cU$ is infinite, and each $x\in X$ belongs to all but finitely many $U\in\cU$.

Let $\cO$, $\Omega$, $\Gamma$ denote the
families of all \emph{open} covers, $\omega$-covers, and $\gamma$-covers of $X$, respectively.
Then $\Gamma\sbst\Omega\sbst \cO$.

For a space $X$ and collections $\scrA,\scrB$ of covers of $X$,
the following properties were introduced by Scheepers in \cite{coc1},
to generalize a variety of classical properties:
\bi
\itm[$\sone(\scrA,\scrB)$:]
For each sequence $\seq{\cU_n}$ of members of $\scrA$,
there exist members $U_n\in\cU_n$, $n\in\N$, such that $\setseq{U_n}\in\scrB$.
\itm[$\sfin(\scrA,\scrB)$:]
For each sequence $\seq{\cU_n}$
of members of $\scrA$, there exist finite
subsets $\cF_n\sbst\cU_n$, $n\in\N$, such that $\Union_{n\in\N}\cF_n\in\scrB$.
\ei
The following notation will also be useful.
\bi
\itm[$\binom{\scrA}{\scrB}$:] Every member of $\scrA$ has a subset which is a member of $\scrB$.
\ei

In 2004, Ko\v{c}inac introduced the properties $\alpha_i$, $i=1,2,3,4$, which are
generalizations of \Arh's $\alpha_i$ local properties.
He initiated their study in \cite{KocAlpha}.
We give the Ko\v{c}inac properties $\alpha_1$ and $\alpha_2$, respectively,
alternative names which are more self-explanatory:
\begin{itemize}
\item[$\scof(\scrA,\scrB)$:]
For each sequence $\seq{\cU_n}$ of members of $\scrA$,
there exist \emph{cofinite} subsets $\cV_n\sbst\cU_n$, $n\in\N$, such that $\Union_{n\in\N}\cV_n\in\scrB$.
\item[$\sinf(\scrA,\scrB)$:]
For each sequence $\seq{\cU_n}$
of members of $\scrA$, there exist \emph{infinite}
subsets $\cV_n\sbst\cU_n$, $n\in\N$, such that $\Union_{n\in\N}\cV_n\in\scrB$.
\end{itemize}

In an independent work \cite{capinf}, we introduced the following selection principle:
\bi
\itm[$\capinf(\scrA,\scrB):$] For each sequence $\seq{\cU_n}$ of elements of $\scrA$,
there is for each $n$ an infinite set $\cV_n\sbst\cU_n$, such that
$\setseq{\bigcap\cV_n}\in\scrB$.
\ei

We classify Ko\v{c}inac's new properties in the case that $\scrA,\scrB$ range over
$\{\cO,\Omega,\Gamma\}$, and describe some relations between them and our
selection principle.
We also classify these properties in the context of topological groups.

Some of the results are stated in a way that makes
them applicable to additional situations.

\section{General topological spaces}

When $\scrA,\scrB$ range over $\{\cO,\Omega,\Gamma\}$,
Ko\v{c}inac's operators $\scof$ and $\sinf$ give
$18$ properties to begin with.

Say that a collection $\cB$ of families of elements of a certain type (e.g., open sets)
is \emph{upward closed} if for each $\cU\in\cB$ and each family $\cV$ of elements of the same
type such that $\cU\sbst\cV$, $\cV\in\cB$.
For example, $\cO,\Omega$ are upward closed, but $\Gamma$ is not.

\begin{lem}\label{upward1}
\mbox{}
\be
\itm For each $\scrB$, $\sinf(\cO,\scrB)$ fails.
\itm If $\scrA\sbst\scrB$ and $\scrB$ is upward closed, then $\scof(\scrA,\scrB)$ holds.
\ee
\end{lem}
\begin{proof}
(1) Each nontrivial (infinite $T_1$) space has a finite cover.

(2) Given elements $\cU_n\in\scrA$, take $\cV_n=\cU_n$ for all $n$.
Clearly, $\Union_n\cV_n\in\scrB$.
\end{proof}

Thus, $\sinf(\cO,\Gamma)$, $\sinf(\cO,\Omega)$, and $\sinf(\cO,\cO)$ always fail,
and $\scof(\Omega,\Omega)$, $\scof(\Gamma,\Omega)$,
$\scof(\Omega,\cO)$, $\scof(\Gamma,\cO)$, and $\scof(\cO,\cO)$ always hold.

\begin{lem}\label{upward2}
Assume that all members of $\scrA$ are infinite,
and $\scrB$ is upwards closed. Then $\scof(\scrA,\scrB)=\sinf(\scrA,\scrB)$.\hfill\qed
\end{lem}

It follows that $\sinf(\Omega,\Omega)$, $\sinf(\Gamma,\Omega)$,
$\sinf(\Omega,\cO)$, and $\sinf(\Gamma,\cO)$ always hold.

The following is immediate.
\begin{lem}\label{choose}
$\sinf(\scrA,\scrB)$ and $\scof(\scrA,\scrB)$, both imply $\binom{\scrA}{\scrB}$.\hfill\qed
\end{lem}

As $\binom{\cO}{\Omega}$ never holds \cite{strongdiags},
$\sinf(\cO,\Omega)$ and $\sinf(\cO,\Gamma)$ always fail.

\begin{lem}
\mbox{}
\be
\itm If $\scrB$ is closed under adding finitely many sets to its elements,
then $\scof(\scrA,\scrB)$ implies $\scrA\sbst\scrB$.
\itm $\scof(\Omega,\Gamma)$ never holds.
\ee
\end{lem}
\begin{proof}
(1) Let $\cU\in\scrA$ and apply $\scof(\scrA,\scrB)$ to the constant sequence
$\cU_n=\cU$, $n\in\N$, to obtain cofinite $\cV_n\sbst\cU$, $n\in\N$, such that
$\cV=\Union_{n\in\N}\cV_n\in\scrB$. $\cV$ is a cofinite subset of $\cU$,
and by the assumption on $\scrB$, $\cU=\cV\cup(\cU\sm\cV)\in\scrB$.

(2) Fix an infinite subset $D$ of $X$. Then $\{X\sm F : \emptyset\neq F\in[D]^{<\alephes}\}\in\Omega\sm\Gamma$.
Use (1).
\end{proof}

Only three properties survive: $\sinf(\Omega,\Gamma)$, $\sinf(\Gamma,\Gamma)$, and
$\scof(\Gamma,\Gamma)$.
These properties are \emph{not} trivial,
as they turn out to characterize known nontrivial properties.

\begin{thm}\label{last3}
\mbox{}
\be
\itm $\sinf(\Gamma,\Gamma)=\sone(\Gamma,\Gamma)$ \cite{KocAlpha}.
\itm $\sinf(\Omega,\Gamma)=\sone(\Omega,\Gamma)$ \cite{KocAlpha}.
\itm Under mild hypotheses on the space $X$, $\scof(\Gamma,\Gamma)=QN$ \cite{BH07, Sakai07}.
\ee
\end{thm}

It follows that no implications are provable among the surviving three properties.

Ko\v{c}inac's proof of Theorem \ref{last3} (1) is
essentially the same as our proof in \cite{capinf} that
$\capinf(\Gamma,\Gamma)=\sone(\Gamma,\Gamma)$.\footnote{Formally, the symbol $\Gamma$
on the right coordinate of $\capinf$ should allow \emph{all} (not necessarily open)
$\gamma$-covers of $X$. We assume that henceforth.}
These two results imply the following, to which we give a direct proof.

\begin{thm}\label{same1}
$\sinf(\Gamma,\Gamma)=\capinf(\Gamma,\Gamma)$.
\end{thm}
\begin{proof}
$(\Impl)$
Assume that $\cU_n$, $n\in\N$, are open $\gamma$-covers of $X$.
As each infinite subset of a $\gamma$-cover is again a $\gamma$-cover,
we may assume that the covers $\cU_n$ are pairwise disjoint.
By $\sinf(\Gamma,\Gamma)$, there are infinite $\cV_n\sbst\cU_n$, $n\in\N$,
such that $\cV=\Union_{n\in\N}\cV_n$ is a $\gamma$-cover of $X$.
Fix $x\in X$. Assume that there are infinitely many $n$ such that
there is $U\in\cV_n$ not containing $x$. As the sets $\cV_n$ are
disjoint, $\cV$ is not a $\gamma$-cover of $X$, a contradiction.
It follows that $\cW=\sseq{\bigcap\cV_n}$ is a $\gamma$-cover of $X$
($\cW$ is infinite because it is an $\omega$-cover of $X$).

$(\Leftarrow)$ Assume that $\cU_n$, $n\in\N$, are open $\gamma$-covers of $X$.
By $\capinf(\Gamma,\Gamma)$, there are infinite $\cV_n\sbst\cU_n$, $n\in\N$,
such that $\sseq{\bigcap\cV_n}$ is a $\gamma$-cover of $X$.
Fix $x\in X$. Let $N$ be such that for all $n\ge N$, $x\in\bigcap\cV_n$.
Then $x\in U$ for all $U\in\cV_n$. Now, for each of the finitely many $n<N$,
$\cV_n$ is an infinite subset of the $\gamma$-cover $\cU_n$ and is
therefore a $\gamma$-cover of $X$. It follows that there are only
finitely many $U\in\cV_n$ such that $x\nin U$. It follows that
$\cV=\Union_{n\in\N}\cV_n$ is a $\gamma$-cover of $X$.
\end{proof}

\begin{cor}\label{happy}
If $\Gamma\sbst\scrA$, then $\sinf(\scrA,\Gamma)=\capinf(\scrA,\Gamma)=\sone(\scrA,\Gamma)$.
\end{cor}
\begin{proof}
Take the conjunction of the properties in Theorem \ref{same1} and the comment
before it with $\binom{\scrA}{\Gamma}$.
\end{proof}

\section{Topological groups}

Ko\v{c}inac also considered in \cite{KocAlpha} the case of topological groups.
Let $G$ be a topological group. For an open neighborhood
$U\neq G$ of the unit $e$, let $o(U)=\{gU : g\in G\}$.
Let $\cO(e)$ be the collection of all these covers $o(U)$.
Ko\v{c}inac asked whether $\sone(\cO(e),\Gamma)$ could hold in
any group. We give a negative answer in a strong sense.

\begin{thm}\label{nogroup}
For each topological group whose topology is nontrivial, $\binom{\cO(e)}{\Omega}$ fails.
\end{thm}
\begin{proof}
Take a neighborhood $U\neq G$ of $e$,
and fix $g\in G\sm U$.
Let $V$ be a neighborhood of $e$ such that $V\cdot V\sbst U$ and $V=V\inv$.

We claim that $o(V)$ is not an $\omega$-cover of $G$. Indeed,
no element of $o(V)$ contains $\{1,g\}$:
If $1\in a\cdot V$, then $a\inv\in V$, hence $a\in V\inv=V$,
and therefore $a\cdot V\sbst V\cdot V\sbst U$. Thus, $g\nin a\cdot V$.
\end{proof}

For $U\sbst G$, let $\omega(U)=\{F\cdot U : F\in [G]^{<\alephes}\}$.
A set $U\sbst G$ is \emph{finitely-bounding} if there is a finite $F\sbst G$
such that $F\cdot U=G$. $\omega(U)$ is an $\w$-cover of $G$ if, and only if,
$U$ is not finitely-bounding.

Let $\Omega(e)$ be the collection of all families $\w(U)$ such that $U$ is an open
neighborhood of $e$ which is not finitely-bounding.
Thus, $\Omega(e)=\emptyset$ if, and only if, $G$ is totally bounded.

We can now classify the group theoretic properties $\Pi(\scrA,\scrB)$ where $\Pi\in\{\sinf,\scof\}$,
$\scrA\in\{\cO(e),\Omega(e)\}$, and $\scrB\in\{\cO,\Omega,\Gamma\}$.

Theorem \ref{nogroup} and Lemma \ref{choose} rule out the $4$ properties
where $(\scrA,\scrB)$ is $(\cO(e),\Omega)$ or $(\cO(e),\Gamma)$.
By Lemma \ref{upward1}, $\scof(\cO(e),\cO)$ always holds.
$\sinf(\cO(e),\cO)$ is also trivial: If there is a finitely bounding
open $U\sbst G$, then $\sinf(\cO(e),\cO)$ fails. And if not,
then $\sinf(\cO(e),\cO)$ holds by Lemmas \ref{upward1} and \ref{upward2}.
By the same Lemmas, $\sinf(\Omega(e),\Omega)$ and $\scof(\Omega(e),\Omega)$
(and therefore also $\scof(\Omega(e),\cO)$ and $\scof(\Omega(e),\cO)$)
always hold.

$\scof(\Omega(e),\Gamma)$ characterizes totally bounded groups.

\begin{thm}\label{tbdd}
The following are equivalent:
\be
\itm $G$ is totally bounded;
\itm $G$ satisfies $\scof(\Omega(e),\Gamma)$.
\ee
\end{thm}
\begin{proof}
$(1\Impl 2)$ If $G$ is totally bounded, then $\Omega(e)=\emptyset$, and therefore $\scof(\Omega(e),\Gamma)$
holds trivially.

$(2\Impl 1)$ Let $U$ be an open neighborhood of $e$ such that for each finite
$F\sbst G$, $F\cdot U\neq G$.
Let $V$ be a neighborhood of $e$ such that $V\cdot V\sbst U$ and $V=V\inv$.
By the proof of Theorem \ref{nogroup}, $o(V)$ is not an $\omega$-cover of $G$.
As $\w(V)$ refines $\w(U)$ and $G\nin\w(U)$, $G\nin\w(V)$ and therefore
$\w(V)$ is an $\w$-cover of $X$, and is therefore infinite. As $\w(V)$ is obtained by taking all
finite unions of elements of $o(V)$, $o(V)$ is infinite.

Assume that $\cV$ is a cofinite subset of $\w(V)$, and that
$\cV$ is a $\gamma$-cover of $G$.
Set $\cW=\cV\cap o(V)$. As $o(V)\sbst\w(V)$ and $\cV$ is cofinite in $\w(V)$,
$\cW$ is cofinite in $o(V)$ and is in particular infinite.
It follows that $\cW$ is an infinite subset of the $\gamma$-cover $\cV$,
and is therefore a $\gamma$-cover of $G$. As $\cW\sbst o(V)$,
$o(V)$ is an $\w$-cover of $G$. A contradiction.
\end{proof}

The only remaining property is $\sinf(\Omega(e),\Gamma)$.
In \cite{KocAlpha} it is shown that $\sinf(\Omega(e),\Gamma)=\sone(\Omega(e),\Gamma)$.
This also follows from Corollary \ref{happy}.
By a result of Babinkostova \cite{Bab05}, for metrizable
groups $\sone(\Omega(e),\Gamma)$ (and therefore $\sinf(\Omega(e),\Gamma)$)
characterizes $\sigma$-totally bounded groups (see also \cite{BG}).
Compare this to Theorem \ref{tbdd}.

\section{Ko\v{c}inac's $\alpha_3$ and $\alpha_4$}

Ko\v{c}inac has also introduced the following properties:
\begin{itemize}
\item[$\alpha_3(\scrA,\scrB)$:]
For each sequence $\seq{\cU_n}$ of members of $\scrA$,
there are an infinite $I\sbst\N$ and \emph{infinite} subsets $\cV_n\sbst\cU_n$, $n\in I$,
such that $\Union_{n\in\N}\cV_n\in\scrB$.
\item[$\alpha_4(\scrA,\scrB)$:]
For each sequence $\seq{\cU_n}$ of members of $\scrA$,
there are an infinite $I\sbst\N$ and \emph{nonempty} subsets $\cV_n\sbst\cU_n$, $n\in I$,
such that $\Union_{n\in\N}\cV_n\in\scrB$.
\end{itemize}

In the two contexts studied here, these selection principles do not give new properties,
but rather are trivial or characterize known properties.

To see this, observe that the following statements hold:
\be
\itm\label{1} If all members of $\scrA$ are infinite, then
$$\sinf(\scrA,\scrB)\Impl\alpha_3(\scrA,\scrB)\Impl\alpha_4(\scrA,\scrB).$$
\itm\label{2} If there are finite elements in $\scrA$, then $\alpha_3(\scrA,\scrB)$ fails.
\itm\label{4} If $\scrA\sbst\scrB$ and $\scrB$ is upward closed, then $\alpha_4(\scrA,\scrB)$ holds.
\itm\label{3} $\alpha_3(\scrA,\scrB),\alpha_4(\scrA,\scrB)\Impl\binom{\scrA}{\scrB}$.
\ee

\eqref{1} trivializes $\alpha_3$ and $\alpha_4$ of the pairs
$(\Omega,\Omega)$, $(\Omega,\Gamma)$, $(\Omega,\cO)$,
and $(\Gamma,\cO)$.
\eqref{2} trivializes $\alpha_3(\cO,\cO)$.
\eqref{4} trivializes $\alpha_4(\cO,\cO)$.
\eqref{3} trivializes $\alpha_3$ and $\alpha_4$ of the pairs
$(\cO,\Omega)$ and $(\cO,\Gamma)$.

$\alpha_3(\Gamma,\Gamma)=\alpha_4(\Gamma,\Gamma)=\sone(\Gamma,\Gamma)$, and
$\alpha_3(\Omega,\Gamma)=\alpha_4(\Omega,\Gamma)=\sone(\Omega,\Gamma)$ \cite{KocAlpha}.

This completes the classification of all mentioned properties for general topological spaces.
The classification of these properties for topological groups is left to the reader.

\ed